\documentclass[twoside]{amsart}
\usepackage{amsmath,amsfonts,amssymb}
\usepackage{mathrsfs}
\usepackage[colorlinks=false,linktocpage=true]{hyperref}
\usepackage{tocvsec2}
\usepackage{latexsym}
\numberwithin{equation}{section}
\font\tenscrpt=eusm10
\font\sevenscrpt=eusm10 scaled 700
\font\fivescrpt=eusm10 scaled 500
\newfam\eusmfam
\textfont\eusmfam=\tenscrpt
\scriptfont\eusmfam=\sevenscrpt
\scriptscriptfont\eusmfam=\fivescrpt
\def\scr#1{{\fam\eusmfam\relax#1}}
\font\eight=cmr8

\def\qed{\quad\vcenter{\hrule\hbox{\vrule height.6em\kern.6em\vrule}\hrule}}

\newenvironment{pf*}[1]{{\textsc #1.\quad}}{$\qed$\bigskip\newline}
\newtheorem{theorem}{Theorem}[section]
\newtheorem{thm}{Theorem}[section]
\newtheorem{cor}{Corollary}[section]
\newtheorem{lem}{Lemma}[section]
\newtheorem{rem}{Remark}[section]
\newtheorem{notn}{Notation}[section]
\newtheorem{defn}{Definition}[section]
\newtheorem{conj}[theorem]{Conjecture}

\newcommand{\thmref}[1]{Theorem~{\rm $\ref{#1}$}}

\newcommand{\lemref}[1]{Lemma~{\rm $\ref{#1}$}}
\newcommand{\corref}[1]{Corollary~{\rm $\ref{#1}$}}

\newcommand{\defnref}[1]{Definition~{\rm $\ref{#1}$}}
\newcommand{\remref}[1]{Remark~{\rm $\ref{#1}$}}
\newcommand{\notnref}[1]{Notation~{\rm $\ref{#1}$}}
\newcommand{\eqnref}[1]{{\rm (\ref{#1})}}

\def\presdh{\mbox{PD}_{B_t,h}}
\def\sdd{\mbox{D}_{B_t,h}}
\def\sdp{{\Bbb D}_{B}}
\def\sds{{\Bbb D}_{B_s}}
\def\sd{{\Bbb D}_{B_t}}

\def\sdS2{{\Bbb D}_{S^{(2)}_t}}
\def\SB{{\scr S}_B(\Rp)}
\def\SBS{{\scr S}_B(\S)}
\def\SBs{{\scr S}_{B}^{\mbox{\eight s}}(\Rp)}

\def\SBSs{{\scr S}_{B}^{\mbox{\eight s}}(\S)}
\def\SBSk{{\scr S}_B^{(k)}(\S)}
\def\SBSks{{\scr S}_B^{\mbox{\eight s},(k)}(\S)}
\def\SBrep{\scr S^{\mbox{\eight{SI}}_B}}
\def\SBcrep{\scr S_c^{\mbox{\eight{SI}}_B}}
\def\SBrepz{\scr S_0^{\mbox{\eight{SI}}_B}}
\def\SBcrepz{{\scr S}_{c,0}^{\mbox{\eight{SI}}_B}}
\def\eqdef{\overset{\triangle}{=}}
\def\f{\frac}
\def\s{\sqrt}
\def\E{\Bbb E}
\def\babs{\begin{abstract}}
\def\eabs{\end{abstract}}
\def\beq{\begin{equation}}
\def\eeq{\end{equation}}
\def\beqs{\begin{equation*}}
\def\eeqs{\end{equation*}}
\def\bsp{\begin{split}}
\def\esp{\end{split}}
\def\bcs{\begin{cases}}
\def\ecs{\end{cases}}
\def\bthm{\begin{thm}}
\def\ethm{\end{thm}}
\def\bcor{\begin{cor}}
\def\ecor{\end{cor}}
\def\bprop{\begin{prop}}
\def\eprop{\end{prop}}
\def\bconj{\begin{conj}}
\def\econj{\end{conj}}
\def\blem{\begin{lem}}
\def\elem{\end{lem}}
\def\bdefn{\begin{defn}}
\def\edefn{\end{defn}}
\def\bexmp{\begin{exmp}}
\def\eexmp{\end{exmp}}
\def\bpf{\begin{proof}}
\def\epf{\end{proof}}

\def\lbr{\left\{}
\def\rbr{\right\}}
\def\l[{\left[}
\def\r]{\right]}
\def\l|{\left|}
\def\r|{\right|}

\def\dis{\displaystyle}

\def\filpsp{(\Omega, \scr{F}, \{{\scr{F}}_t\},\P)}
\def\Ft{\scr{F}_t}

\def\eqdef{\overset{\triangle}{=}}

\def\P{{\mathbb P}}

\def\EP{{\mathbb E}_{\P}}

\def\R{{\mathbb R}}

\def\S{{\mathbb S}}
\def\Rp{{\mathbb R}_+}
\def\Rpo{{\mathbb R}_+^0}

\def\Rpo{{\overset{\circ}{{\mathbb R}}_{+}}}

\def\Filt{{\scr F}_t}
\def\P{{\mathbb P}}

\def\EP{{\mathbb E}_{\P}}

\def\BRp{{\scr B}(\Rp)}

\def\R{{\mathbb R}}
\def\S{{\mathbb S}}

\begin{document}
\title{A Differentiation Theory for It\^o's Calculus}\footnote{This is a preprint of an article published in the [Stochastic Analysis and Applications] [2006] [copyright Taylor \& Francis]; available online at:\\  http://www.informaworld.com/smpp/6561327-36708330/content~db=all~content=a742126397}
\author{Hassan Allouba}
\date{7/20/2004}
\keywords{It\^o calculus; SDEs; SPDEs; Stochastic calculus; Stochastic derivative.}
\subjclass[2000]{60H05; 60H10; 60H15; 60G20; 60G05}
\maketitle
\begin{abstract}
A peculiar feature of It\^o's calculus is that it is an integral calculus that gives
  no explicit derivative with a systematic differentiation theory counterpart, as in elementary
  calculus.   So,  can we define a pathwise stochastic derivative of semimartingales with respect to
  Brownian motion that leads to a differentiation theory counterpart to It\^o's integral calculus?
  From It\^o's definition of his integral, such a derivative must be based on the quadratic covariation process.
  We give such a derivative in this note and we show that it leads to a fundamental theorem of stochastic calculus,
   a generalized stochastic chain rule that includes the case of convex functions acting on continuous semimartingales,
   and the stochastic mean value and Rolle's theorems.   In addition, it interacts with
  basic algebraic operations on semimartingales similarly to the way the deterministic derivative
  does on deterministic functions, making it natural for computations.   Such a differentiation theory leads
  to many interesting applications some of which we address in an upcoming article.
\end{abstract}

  \section{Differentiating Semimartingales With Respect to Brownian Motion}
  One of the greatest  twentieth century's discoveries in probability and mathematics
  is It\^o's theory of stochastic integration \cite{Ito} which, in its simplest form, shows how to integrate
  certain stochastic processes with respect to Brownian motion (BM).  It\^o's powerful ideas are
  still at the heart of some of the most important advances in both pure and applied mathematics sixty years later
  (stochastic analysis, SDEs,  SPDEs, finance, and others).    As is well known, It\^o's calculus is an integral calculus
  that gives no explicit derivative and no systematic differentiation theory counterpart, as in elementary
  calculus.   So, the question is:  can we define a pathwise stochastic derivative of semimartingales with respect to
  BM that leads to a differentiation theory counterpart to It\^o's integral calculus?
  Before giving such a derivative, we briefly recall that the essential ingredient in It\^o's definition
  of the integral $\int_0^tX_sdB_s$ of a stochastic process $X=\{X_t;0\le t<\infty\}$ with respect to a
  BM $B=\{B_t;0\le t<\infty\}$ is It\^o's isometry
  $$\E\left(\int_0^tX_s dB_s\right)^2=\E\int_0^tX_s^2ds,$$
  where $\E$ is the usual mathematical expectation.  This isometry leads to a definition of the integral with respect
  to $B$ in terms of one with respect to the
  quadratic variation of $B$ given by $\langle B,B\rangle_t(\omega)=\langle B\rangle_t(\omega)=t$.  Our idea is then to
  define the derivative $dS_t/dB_t$ of a semimartingale $S$ with respect to a BM $B$ at $t$ as a generalized version of
  the pathwise stochastic derivative
   $d\langle S,B\rangle_{t}(\omega)/d\langle B\rangle_{t}(\omega)=d\langle S,B\rangle_{t}(\omega)/dt$
   of the covariation of $S$ and $B$ at $t$
   with respect to the quadratic variation of $B$ at $t$.   In fact,
  our derivative covers cases where  $d\langle S,B\rangle_{t}(\omega)/dt$
  doesn't exist for all $t$.   We use an integral formulation to define our stochastic derivative (\defnref{Ad}), which allows
  for greater applicability: for example if $t\mapsto\langle S,B\rangle_t$ is convex almost surely, then
  our derivative of $S$ with respect to Brownian motion $B$ exists for all $t\in\Rp$ a.s.~(see \lemref{1}).
  We show that our derivative is an anti It\^o integral that has the ``desired'' properties leading to
   a fundamental theorem of stochastic calculus (\thmref{ft}), a generalized stochastic chain rule
   (\thmref{cr}) that includes the case of convex functions acting on continuous semimartingales,
   and the stochastic mean value and Rolle's theorems (\lemref{1}).   In addition, it interacts with
  basic algebraic operations on semimartingales similarly to the way Newton's deterministic derivative
  does on deterministic functions (\thmref{rules}), making it natural for computations.

  Heuristically, our derivative of a semimartingale $S$ with respect to a BM $B$ is the pathwise
  ``velocity'' of $S$ relative to $B$ at each point in time (we call it the $B$-Brownian velocity of $S$, see
  the examples at the end of Section 3).
  It gives a Brownian path view of the changes in a
  semimartingale's path $S(\omega)$ in time by measuring the rate of change of the covariation of $S$
  with the Brownian motion $B$ with respect to the quadratic variation of $B$ (time).
  The sign of our derivative of $S$ with respect to $B$ tells us, path by path, whether $S$ and $B$ are increasing and
  decreasing together (positive sign) or whether the value of $S$ changes in the opposite direction of changes in the value
  of $B$ (negative sign).   The magnitude of the $B$-Brownian velocity of $S$ gives us the $B$-Brownian speed of $S$.

  We believe that this pathwise view, in addition to giving rise to a differentiation theory for It\^o's calculus,
   is a useful tool in the analysis of SDEs and SPDEs; leading to a new
  smoothness and regularity theory for solutions of these stochastic equations, which
  in turns leads to new insights into the behavior of
  solutions to SDEs and SPDEs relative to their driving noise.  It also leads to
  a rich differential stochastic vector calculus
  as well as to a new stochastic optimization theory (optimization with respect to the noise).      We
  study these different aspects in more details in \cite{A} and in subsequent articles.
  \begin{notn}   \label{notn} Throughout this article,  let $S=\{S_t,\Ft;t\in\Rp\}$ be
  a continuous semimartingale, let $M=\{M_t,\Ft;t\in\Rp\}$ and $V=\{V_t,\Ft;t\in\Rp\}$
  be the continuous local martingale and the continuous process of bounded
  variation in the decomposition of  $S$, respectively, and let $B=\{B_t,\Ft;t\in\Rp\}$ be
  a standard BM on the same usual probability space $\filpsp$ (a usual probability space is
  one where the filtration $\{\Filt\}$ satisfies the usual conditions: right continuity and completeness).
  Also, let $D\langle S,B\rangle_{t}\eqdef d\langle S,B\rangle_{t}/dt$, which we call the strong derivative
  of $S$ with respect to $B$.
  Finally, we denote by $\lambda$ the Lebesgue measure
  on $\Rp$ and by $\Rpo$ the set $\Rp\backslash\{0\}$.
  \end{notn}

  \begin{defn}[Derivative of semimartingale with respect to BM]\label{Ad}
    The \textit{stochastic difference} and \textit{stochastic derivative} of $S$ with
    respect to $B$ at $t$, respectively, are defined by
    \beq
    \begin{split}
    \sdd S_t&\eqdef\begin{cases}
    \dis{\f{3}{2h^3}}\int_0^h r\left[\langle S,B\rangle_{t+r}
    -\langle S,B\rangle_{t-r}\right]dr;&0<t<\infty,\ h>0,\\ \vspace{-4.5mm}  \\
     \dis{\f{3}{h^3}}\int_0^h r\langle S,B\rangle_{r}\, dr;&t=0,\ h>0.
    \end{cases}
    \\
    \sd S_t&\eqdef\lim_{h\to0} \sdd S_t; \ 0\le t<\infty,
    \end{split}
    \label{sd}
    \eeq
    whenever this limit exists.  We call continuous
    semimartingales for which the limit in \eqref{sd}
    exists for all $t\in\S\in\BRp$ a.s.~differentiable with respect to $B$, on $\S$,
    and we denote this class by $\SBS$.
    The $k$-th $B$-derivative of $S$ is defined iteratively in the
    obvious way, and the class of $k$-times $B$-differentiable elements of
    $\SBS$ is denoted by $\SBSk$.   If the derivative $d\langle S,B\rangle_t/dt$
    exists then
    \begin{equation}
    \sd S_t= \frac{d\langle S,B\rangle_t}{dt}
    \label{sB}
    \end{equation}
    $($see \thmref{1} below$)$,
    and we call $d\langle S,B\rangle_t/dt$ the strong derivative of $S$ with respect
    to $B$ at $t$.  We denote the class of continuous semimartingales whose strong
    $B$-derivative exists on $\S$ by $\SBSs$.  The class of $k$-times strongly
    $B$-differentiable elements of $\SBSs$ is denoted by $\SBSks$.

    If the strong $B$-derivative of $S$ exists at $t$ and if $f\in C^1(\R;\R)$ with $f'(x)\neq0$
    for all $x\in\R$ and the map $x\mapsto f'(x)$ is absolutely continuous,
    we define the strong derivative of $S$ with respect to the semimartingale $S^{(2)}\eqdef f(B)$ at $t$ by
    \begin{equation}
    \sdS2 S_t \eqdef \frac{d \langle S,B\rangle_t}{d \langle S^{(2)}\rangle_t}.
    \label{sfB}
    \end{equation}
    We also have the same definition of $\sdS2 S_t$ in the case $f$ is convex and $f'_-(x)\neq0$
    for all $x\in\R$, where   $f'_-(x)$ is the left derivative of $f$ at $x$.
    The generalized derivative of $S$ with respect to $f(B)$ at $t$ is obtained straightforwardly from \eqref{sd}.
        \end{defn}
     \begin{rem}\label{remsd}
    $(a)$ It follows immediately from our \defnref{Ad} of the stochastic derivative process that $\sdp V\equiv0$
    for all processes $V$ of bounded variation on compacts: there changes are ``too slow'' for the BM to pickup,
    and they behave like constants in elementary calculus in that their derivative is $0$.    Additionally,
    if $M$ and $B$ are independent or orthogonal $(\langle M,B\rangle\equiv0)$;
    then from the definitions above  $\sdp M\equiv0$$:$
    $M$ is independent of or orthogonal to $B$ and therefore it is ``unaffected'' by changes in $B$.  \\ \vspace{-3mm}  \\
    $(b)$ We don't use the following observation here in this note; but, we formally
    think of $\sdd S_t$ as the     quadratic variation of a generalized stochastic integral which we call
    the pre-stochastic difference of $S$ with respect to the Brownian motion
     $B$:
    \begin{equation*}
     \presdh S_t\eqdef \begin{cases}\dis\s{\f{3}{2h^3}}\int_0^h\s r\left[\langle S,B\rangle_{t+r}
    -\langle S,B\rangle_{t-r}\right]^{\f12}dB_r;&0<t<\infty,\ h>0,\\
     \dis\s{\f{3}{h^3}}\int_0^h\s r\left[\langle S,B\rangle_{r}
   \right]^{\f12}dB_r;&t=0,\ h>0.
       \end{cases}
    \label{psd}
    \end{equation*}
        \end{rem}

    \section{Fundamental Results}
    We start with a lemma which expresses our stochastic derivative as
    a generalized derivative of the covariance process $\langle M,B\rangle$ with respect to
    time $t=\langle B\rangle_t$ as well as gives us
    a stochastic mean value theorem (SMVT) and a stochastic Rolle's theorem (SRT).
    \begin{lem}\label{1}
     Let $S$, $M$, $V$, and $B$ be defined as in \notnref{notn}
     on the same usual probability space $\filpsp$.
     \begin{enumerate}\renewcommand{\labelenumi}{$(\alph{enumi})$}
     \item $($Stochastic derivative as a generalized derivative of the covariance process$)$
      Assume that the one sided derivatives processes $D^+\langle M,B\rangle$ and
      $D^-\langle M,B\rangle$ are finite a.s.;  i.e.,
    \begin{equation}
    \begin{cases}
    \dis-\infty<D^\pm\langle M,B\rangle_{t}
    =\lim_{\epsilon\to0^\pm}\f{\langle M,B\rangle_{t+\epsilon}-\langle M,B\rangle_{t}}{\epsilon}<\infty;&
    0< t<\infty\mbox{ a.s.~}\P,\\
     -\infty<\left(D^+\langle M,B\rangle_{t}\right)|_{t=0}<\infty;&t=0\mbox{ a.s.~}\P
     \end{cases}
     \label{fosd}
     \end{equation}
   $($e.g., if $t\mapsto\langle M,B\rangle_t$ is convex a.s.~$\P).$
    Then, $S\in\SB$ and the stochastic derivative of $S$ with respect to $B$
    is given by
    \begin{equation}
    \begin{split}
     \sd S_t=\begin{cases}
    \dis \f12\left[D^+\langle M,B\rangle_{t}
    +D^-\langle M,B\rangle_{t}\right];& \ 0<t<\infty, \mbox{ a.s.~}\P,\\
     \left(D^+\langle M,B\rangle_{t}\right)|_{t=0};&t=0,\mbox{ a.s.~}\P.
     \end{cases}
    \end{split}
    \label{davg}
    \end{equation}
    In particular, if \eqref{davg} holds on $\Omega^*\subset\Omega$ $($with $\P(\Omega^*)=1$$)$
    and $\langle M,B\rangle(\omega_0)$ is differentiable at $t$ for some $\omega_0\in\Omega^*$; then
    $\sd S_t(\omega_0)=D\langle M,B\rangle_t(\omega_0)$.\\  \\
    \item $($SMVT and SRT$)$  Let $\langle M,B\rangle$ be continuous on the
    closed interval $[a,b]\subset\Rp$ and differentiable on the open interval $(a,b)$, a.s.~$\P;$
    then,
     \beq
     \left(\sd S_t\right){\mid_{t=c}}(\omega)=\f{\langle M,B\rangle_{b}(\omega)
     -\langle M,B\rangle_{a}(\omega)}{b-a}, 
     \label{smvt}
     \end{equation}
     for some random variable $c(\omega)\in(a,b)$ a.s. $\P$.
     In particular, if $\langle M,B\rangle_{b}(\omega)=\langle M,B\rangle_{a}(\omega)$
     a.s.~$\P;$ then $\sd S_t{\mid_{t=c}}(\omega)=0$ a.s.~$\P$.
     \end{enumerate}
     \end{lem}
     A simple application of the SMVT leads to
     \begin{cor}  \label{cons}
     Suppose that, for some BM $B$, $S\in\SBs$ with decomposition $S_t=S_0+V_t+M_t$, $t\ge0$, a.s.~$\P$.
     \begin{enumerate}\renewcommand{\labelenumi}{(\roman{enumi})}
     \item
    If $\sd S_t=0$ for all $t>0$ a.s. $\P$.  Then,
     $\langle M,B\rangle\equiv0$ a.s.~$\P$.
     \item If $\sd S_t$ does not change sign on $(a,b)$ a.s.~$\P$, then
     $\langle M,B\rangle$ is monotonic over $(a,b)$ a.s.~$\P$: increasing, nondecreasing,
     decreasing, or nonincreasing as
     $$\sd S_t>0,\ \sd S_t\ge0,\ \sd S_t<0,\ \mbox{or}\ \sd S_t\le0; \forall\ t\in(a,b)\ \mbox{a.s.}~\P,$$
     respectively.
     \item If $$\left|\sd S_t\right|\le K;\quad a<t<b, \mbox{ a.s. }\P$$
     then $\langle M,B\rangle$ is Lipschitz on $(a,b):$ 
     $$\left|\langle M,B\rangle_t-\langle M,B\rangle_s\right|\le K|t-s|;\quad t,s\in(a,b)\mbox{ a.s. }\P.$$
     \end{enumerate}   
     \end{cor}
      We state the Fundamental theorem of stochastic calculus for a class of semimartingales that covers
    It\^o SDEs of interest and that can be generalized to cover SPDEs of interest as well:
   \begin{equation}
   S=\left\{S_t=S_0+V_t+\int_0^tX_sdB_s,\Ft;\ 0\le t<\infty\right\}
   \label{sem}
   \end{equation}
   where $V$ and $B$ are as in \notnref{notn};
   and  the adapted process $X=\{X_t,\Ft;t\in\Rp\}$ is in $L^2([0,t];\lambda)$ $\forall t>0$ a.s.~$\P$.
   We denote  by $\SBrep$ the class of continuous semimartingales $S$ whose local martingale part $M$
   is given by the stochastic integral in \eqnref{sem}; i.e.,
    \begin{equation}
   M_t= \int_0^tX_sdB_s,\Ft;\ 0\le t<\infty.
   \label{rep}
   \end{equation}
   The elements of $\SBrep$ in which the integrand $X$ has a.s.~continuous paths form the
   subclass which we denote by $\SBcrep$.  Finally,
   we denote by $\SBrepz$ and $\SBcrepz$ the subclasses $\SBrepz\subset\SBrep$ and
   $\SBcrepz\subset\SBcrep$ in which $V\equiv0$ for all of its elements.

    We are now ready to present the pathwise fundamental theorem of stochastic calculus (FTSC)
    \begin{theorem}[FTSC]  \label{ft} Let $S\in\SBcrep$.    Then, $S\in\SBs$---in particular
     the process $\langle M,B\rangle$ is differentiable for all $t\in\Rp$, a.s.~$\P$---and
    \begin{enumerate}\renewcommand{\labelenumi}{$(\roman{enumi})$}
    \item the stochastic derivative process
    $\sdp S=\lbr\sd S_t,\Ft;t\in\Rp\rbr$ is given by $ \sd S_t=X_t\mbox{ for all } t\in\Rp,$ a.s.~$\P$.
    In particular,
    \begin{equation}
    \sd\int_0^tX_sdB_s=X_t;\ \forall\,0\le t<\infty,\ \mbox{a.s.}~\P.
    \label{F1s}
    \end{equation}
    Moreover, if $S,\tilde{S}\in\SBcrep$ with $\P\lbr M_t=\tilde{M}_t;\ \forall\,t\in\Rp\rbr=1;$
    then their stochastic derivative processes are indistinguishable$:\ \P\lbr\sd S_t=\sd \tilde S_t;\ \forall\,t\in\Rp\rbr=1$.
     \item
    \begin{equation}
    \int_0^t \sds S_sdB_s=\tilde S_t-\tilde S_0-\tilde V_t; \ \forall\,0\le  t<\infty,\ \mbox{a.s.}~\P
    \label{FIIs}
    \end{equation}
    for any $\tilde{S}\in\SBcrep$ whose local martingale part $\tilde M$ is indistinguishable
    from $M$.  In particular; $\int_0^t \sds S_sdB_s=S_t- S_0- V_t \ \mbox{for all } t\in\Rp$,
     a.s.~$\P$.
    Thus, if $S\in\SBcrepz; \mbox{ then, }$ $\int_0^t \sds S_sdB_s=S_t-S_0\ \mbox{for all } t\in\Rp$,
    a.s.~$\P$.
    \end{enumerate}
     \end{theorem}
      \begin{rem}\label{2}
      In contrast to the fundamental theorem of deterministic calculus, part (ii) of
      \thmref{ft} involves the additional term of bounded variation on compacts
      $\tilde V$, unless $S\in\SBcrepz$.  Remember
      however that, a.s.~$\P$, $\sdp \tilde V\equiv0$ by \remref{remsd}.
       \end{rem}

       In the case the integrand $X$ is not necessarily continuous we state the following
       \begin{theorem} \label{bmart}
        Assume that $S\in\SBrep$; then, a.s.~$\P$,
        the process $\langle M,B\rangle$ is differentiable for all $t\in\Rp\backslash Z$
        and $\sdp S=\lbr X_t;t\in\Rpo\backslash Z\rbr$, for some $Z(\omega)\subset\Rp$
        with  $\lambda(Z)=0$.  If, additionally, the condition \eqref{fosd} hold; then we also have
\begin{equation}
    \begin{split}
     \sd S_t=\begin{cases}
    \dis \f12\left[D^+\int_0^{t}X_s\,ds+D^-\int_0^{t}X_s\,ds\right];&
    \ t\in(0,\infty)\cap Z,\\
     \dis\left(D^+\int_0^{t}X_s\,ds\right)|_{t=0};&t=0,
     \end{cases}
    \end{split}
    \label{davgN}
    \end{equation}
    a.s.~$\P$.
    If $S,\tilde S\in\SBrep$ with indistinguishable $M$ and $\tilde M$;
        then, a.s.~$\P$,
        $\sd S_t=\sd \tilde S_t \mbox{ for all } t\in\Rpo\backslash O$,
        for some $O(\omega)\subset\Rp$ with $\lambda(O)=0$ $($we say that 
        $\sdp S$ and $\sdp \tilde S$ are almost indistinguishable$)$.
        In particular, if $M$ is a $B$-Brownian martingale; then, a.s.~$\P$,
        $\sdp M=\{Y_t;t\in\Rpo\backslash Z\}$
       for some $Z(\omega)\subset\Rp$
        with  $\lambda(Z)=0$, where $Y$ is the progressively measurable process such that $\EP\int_0^tY_s^2ds<\infty$ and
       $\int_0^t Y_s dB_s=M_t$, for all $t\in\Rp$.
             \end{theorem}
     \section{Pathwise Derivative Rules}
     We now show that our pathwise derivative for It\^o's calculus generalizes familiar differentiation rules from
     deterministic to stochastic calculus, making it useful for computations involving functions of semimartingales
     and algebraic operations on several semimartingales.
     We start with our chain rule for stochastic calculus.
      \begin{theorem}[The chain rule of stochastic calculus]  \label{cr}
       \begin{enumerate}\renewcommand{\labelenumi}{$(\alph{enumi})$}
     \item Suppose that $f\in C^1(\R;\R)$ such that the function $x\mapsto f'(x)$ is absolutely continuous.
    \begin{enumerate}\renewcommand{\labelenumii}{$(\roman{enumii})$}
    \item Then, a.s.~$\P$, the process $\sdp f(S)=\{\sd f(S_t);t\in\Rp\}$ is given by
       \begin{equation}
       \sd f(S_t)=\begin{cases}
       \dis\f12\left[\left(D^+\int_0^tf'(S_s)\,d\langle M,B\rangle_{s}\right)
    +\left(D^-\int_0^tf'(S_s)\,d\langle M,B\rangle_{s}\right)\right]; \ 0<t<\infty,\\
     \dis\left.\left( D^+\int_0^tf'(S_s)\,d\langle M,B\rangle_{s}\right)\right|_{t=0};\ t=0,
     \end{cases}
       \label{GCR}
       \end{equation}
       whenever the one sided derivatives are finite.
       If $d\langle M,B\rangle_s=X_{M,B}(s)ds$ and $X_{M,B}$
       has continuous paths on $\Rp$ a.s.~$\P$,
        then $S\in\SBs$ and \eqref{GCR} becomes
       \begin{equation}
       \sd f(S_t)=f'(S_t)\sd S_t=f'(S_t)X_{M,B}(t); \ 0\le t<\infty \mbox{ a.s.}~\P.
       \label{GCR2}
       \end{equation}
       In particular, if $S\in\SBcrep$; then
       \begin{equation}
       \sd f(S_t)=f'(S_t) \sd S_t; \ 0\le t<\infty \mbox{ a.s.}~\P.
       \label{CR}
       \end{equation}
       \item  Suppose further that $f'(x)\neq0$ for every $x\in\R$, $d\langle M,B\rangle_s=X_{M,B}(s)ds$ and $X_{M,B}$
       has continuous paths on $\Rp$ a.s.~$\P$, and let $S^{(2)}\eqdef f(B)$.  Then, $S,S^{(2)}\in\SBs$ and
        \begin{equation}
    \sd S_t = \sdS2 S_t\cdot\sd S^{(2)}(t)\ 0\le t<\infty \mbox{ a.s.}~\P.
    \label{cr2}
    \end{equation}
       \end{enumerate}
       \item If $f:\R\to\R$ is only assumed to be convex;  then \eqref{GCR}, \eqref{GCR2}, and \eqref{CR}
       all hold, replacing $f'(x)$ by the left derivative at $x$, $f'_-(x)$.   Moreover, if in addition the assumptions in part $(a) (ii)$
       hold $($again replacing $f'(x)$ by $f'_-(x))$, then \eqref{cr2} holds.
       \end{enumerate}
            \end{theorem}

       An interesting question then is when is $\sdp S$ itself a martingale (or a local martingale)? It is for example
       clear from the above discussion that $\{\sd (B_t^2-t)=2B_t,\Ft; t\ge0\}$ is a martingale.  The next corollary, which follows as an immediate
       consequence of \thmref{cr}, gives a sufficient condition.
       \begin{cor} \label{derc2}
       Suppose $f\in C^1(\R;\R)$ such that the function $x\mapsto f'(x)$ is absolutely continuous
       and $S_t=f(B_t) +V_t$, $t\ge0$; where $B$ and $V$ are as in \notnref{notn}.   Then
       $$\sd S_t=f'(B_t); \quad0\le t<\infty.$$
       In particular, $\sd S=\{\sd S_t,\Ft;0\le t<\infty\}$ is a martingale $($local martingale$)$ iff the process
       $\{f'(B_t),\Ft;0\le t<\infty\}$ is a martingale $($local martingale$)$.
       If $f:\R\to\R$ is only assumed to be convex, then
       $\sd S=\{\sd S_t,\Ft;0\le t<\infty\}$ is a martingale $($local martingale$)$ iff
       $\{f'_-(B_t),\Ft;0\le t<\infty\}$ is a martingale $($local martingale$)$.
             \end{cor}
       Another immediate consequence of \thmref{cr} is the power rule for our
       pathwise derivative
       \begin{cor}[The power rule of stochastic calculus]  \label{power}
        Let $S\in\SB$ with continuous local martingale part $M$.
        If $d\langle M,B\rangle_s=X_{M,B}(s)ds$ and $X_{M,B}$ has continuous paths on $\Rp$  a.s.~$\P$,
        then $S\in\SBs$ and
        \begin{equation}
        \sd\left(S_t\right)^p=p\left(S_t\right)^{p-1}\sd S_t; \quad 0\le t<\infty,\ \forall\ p\ge1.
        \label{pow1}
         \end{equation}
        If in addition $S_t\neq0$ for every $t\in\Rp$, a.s.~$\P$, then
         \begin{equation}
        \sd\left(S_t\right)^p=p\left(S_t\right)^{p-1}\sd S_t; \quad0\le t<\infty,\ \forall\ p\in\R.
        \label{pow2}
         \end{equation}
            \end{cor}

\begin{theorem}[The sum, product, and ratio rules ] \label{rules} Let $S_1,S_2\in\SB$ and
let $a,b\in\R$ be arbitrary but fixed, then $a S^{(1)}\pm b S^{(2)}\in\SB$ and
\begin{equation}
\begin{split}
\sd \left(a S^{(1)}_t\pm b S^{(2)}_t\right)=a \sd S^{(1)}_t\pm b \sd S^{(2)}_t;\quad0\le t<\infty, \mbox{ a.s.}~\P.
\end{split}
\label{as}
\end{equation}
If in addition the continuous local martingale parts $M^{(1)}$ and $M^{(2)}$ of $S^{(1)}$ and $S^{(2)}$, respectively,
satisfy
\begin{equation}
d\langle M^{(i)},B\rangle_s=X_{M^{(i)},B}(s)ds; \quad i=1,2,
\label{cvcond}
\end{equation}
and $X_{M^{(i)},B}$ has continuous paths on $\Rp$ a.s.~$\P$ for $i=1,2$ then $S^{(1)},S^{(2)},S^{(1)}S^{(2)}\in\SBs$ and
\begin{equation}
\begin{split}
\sd\left (S^{(1)}_t S^{(2)}_t\right)=S^{(2)}_t\sd S^{(1)}_t+ S^{(1)}_t\sd S^{(2)}_t;\quad0\le t<\infty, \mbox{ a.s.}~\P.
\end{split}
\label{m}
\end{equation}
If in addition $S^{(2)}_t\neq0$ for every $t\in\Rp$, a.s.~$\P$, then $S^{(1)}/S^{(2)}\in\SBs$ and
\begin{equation}
\begin{split}
\sd \dis\left(\frac{ S^{(1)}_t}{ S^{(2)}_t}\right)=\frac{S^{(2)}_t\sd S^{(1)}_t- S^{(1)}_t\sd S^{(2)}_t}{\left[S^{(2)}_t\right]^2};
\quad0\le t<\infty, \mbox{ a.s.}~\P.
\end{split}
\label{d}
\end{equation}
\end{theorem}

We now briefly look at the Brownian derivatives (Brownian velocity)
of several simple semimartingales obtained as simple applications
of \thmref{ft}, \thmref{cr}, and \thmref{rules}.  These examples show that
our derivative gives intuitive answers:
      \begin{enumerate}
      \item $$\sd\left|B_t\right|=\begin{cases}1, &B_t>0\\-1&B_t\le0 \end{cases}$$
      I.e., the Brownian speed of $\left|B_t\right|$ is $1$ for all $t\in\Rp$, and the direction of change of $\left|B\right|$
      (increase or decrease) relative to $B$ at $t$ is the same as $B_t$ if $B_t>0$ and is opposite to that of $B_t$ if $B_t<0$.
      \item Let $S_t=B_t^2-V_t$, for any process $V$ as in \notnref{notn}, then $\sd S_t=2B_t$ for all $t\in\Rp$.
      So that the Brownian speed of $S_t$ at any $t\in\Rp$ is $2|B_t|$, and the direction of change of $S_t$
      (increase or decrease) relative to $B$ at $t$ is the same as $B_t$ if $B_t>0$ and is opposite to that of $B_t$ if $B_t<0$.
      \item Let $$\Xi_t^{X,B}=\exp{\left(\int_0^tX_sdB_s-\frac12\int_0^tX^2_sds\right)},$$
       and let the adapted process $X=\{X_t,\Ft;t\in\Rp\}$ be continuous (thus in $L^2([0,t];\lambda)$ $\forall t>0$ a.s.~$\P$).
      Then, $\sd\Xi_t^{X,B}=X_t\Xi_t^{X,B}$.   So that the Brownian speed of the exponential local martingale $\Xi_t^{X,B}$
      at any $t\in\Rp$ is $|X_t|\Xi_t^{X,B}$, and  the direction of change of $\Xi_t^{X,B}$
      (increase or decrease) relative to $B$ at $t$ is the same as $B_t$ if $X_t>0$ and is opposite to that of $B_t$ if $X_t<0$.
      \end{enumerate}
      More applications of our theory presented here (including SDEs and SPDEs) are dealt with in \cite{A} and subsequent articles.
\section{Proofs of results}
\begin{pf*}{Proof of \lemref{1}}
  \begin{enumerate}\renewcommand{\labelenumi}{$(\alph{enumi})$}
     \item  Under the assumptions given, the conclusion follows
from the definition of $\sd S_t$, the facts that $\langle S,B\rangle_{t}=\langle M,B\rangle_{t}$
and $\langle M,B\rangle_{0}=0$,
and the fact that the generalized derivative
\begin{equation}
\scr D g(t)\eqdef
\begin{cases}
\dis{\lim_{h\to0}\f{3}{2h^3}}\int_0^h r\left[g(t+r)-g(t-r)\right]dr
=\f12\left[D^+g(t)+D^-g(t)\right];&0<t<\infty,\\ \vspace{-4.5mm}\\
     \dis{\lim_{h\to0}\f{3}{h^3}}\int_0^h r\left[g(r)-g(0)\right]dr=D^+g(0);&t=0.
\end{cases}
\label{gdir}
\end{equation}
for any function $g:\Rp\to\R$ with finite one sided derivatives $D^\pm g(t)$ for $0<t<\infty$ and
with finite right hand derivative at zero $D^+g(0)$.  Of course if $g$ is differentiable at $t$ then
$\scr D g(t)=g'(t)$.\\ \vspace{-2mm}\\
\item Let the assumptions hold on $\Omega^*\subset\Omega$, with $\P(\Omega^*)=1$, and
fix $\omega\in\Omega^*$.  Then, by the classical mean value theorem $\exists$ a $c(\omega)\in(a,b)$
$\ni$
$$D\langle M,B\rangle_{t}{\mid_{t=c}}(\omega)=\f{\langle M,B\rangle_{b}(\omega)
     -\langle M,B\rangle_{a}(\omega)}{b-a}.$$
Since $\langle M,B\rangle$ is assumed differentiable on $(a,b)$ $\forall\omega\in\Omega^*$; then, by part (a),
$ \sd S_t{\mid_{t=c}}(\omega)=D\langle M,B\rangle_{t}{\mid_{t=c}}(\omega)$ for all $\omega\in\Omega^*$
and \eqnref{smvt} is established.
\end{enumerate}
The proof is complete.\end{pf*}

\begin{pf*}{Proof of \thmref{ft}}
$(i)$
    Under the assumptions on $S$, the stochastic difference and stochastic derivative of $S$ with
    respect to $B$ are  a.s.~$\P$ given, respectively, by
    \beq
    \begin{split}
    \sdd S_t&={\f{3}{2h^3}}\int_0^h r\left[\int_0^{t+r}X_sds
    -\int_0^{t-r}X_sds\right]dr; \quad0<t<\infty,\\
     \sdd S_t&={\f{3}{2h^3}}\int_0^h r\left[\int_0^{r}X_sds\right]dr; \quad t=0,\\
    \sd S_t&=\lim_{h\to0} \sdd S_t=X_t;\quad 0\le t<\infty.
    \end{split}
    \label{sdW}
    \eeq
    where the last equality follows by letting $g(t)=\int_0^{t}X_sds$ and using the
    fundamental theorem of classical calculus along with continuity of $X$ and
    the fact \eqnref{gdir}.
    Now, if $S,\tilde{S}\in\SBcrep$ with $M\equiv\tilde{M}$ a.s.~$\P$, then
    $\int_0^t(X_s-\tilde X_s)^2ds=0\ \forall\,t\in\Rp$ a.s.~$\P$; which, by the continuity
    of $X$ and $\tilde X$ and the first part of the proof, implies that $\sdp S\equiv X\equiv \tilde X
    \equiv\sdp \tilde S$ a.s.~$\P$.
    \\ \vspace{-2mm}\\
$(ii)$
    By part $(i)$ and \eqref{sem} and \eqref{rep} we have $\sdp S\equiv X$ a.s.~$\P$,
    from which \eqref{FIIs} follows.  The rest of the assertions follow immediately,
    completing the proof.
\end{pf*}
\begin{pf*}{Proof of \thmref{bmart}}
   The proof proceeds exactly as in the proof of part (i) of \thmref{ft},
   taking into account the noncontinuity of $X$ and using the fundamental theorem of Lebesgue calculus
   (e.g., see Theorem 10 on p.~107 in \cite{R}) and \eqref{davg} in \lemref{1}.
\end{pf*}
 \begin{pf*}{Proof of \thmref{cr}}    \begin{enumerate}\renewcommand{\labelenumi}{$(\alph{enumi})$}
     \item If  the function $x\mapsto f'(x)$ is absolutely continuous, then
 $f''$ exists Lebesgue-almost everywhere and we have that the It\^o formula
 for $f(S_t)$ is given by
 \begin{equation}
f(S_t)=f(S_0)+\int_0^tf'(S_s)[dM_s+dV_s]+\f12\int_0^tf''(S_s)d\langle M \rangle_s;\quad t\in\Rp,\mbox{ a.s.}~\P.
 \label{Ito}
 \end{equation}
 \begin{enumerate}\renewcommand{\labelenumii}{$(\roman{enumii})$}
    \item Using It\^o's rule \eqref{Ito}, \remref{remsd} (a) ($\sd U\equiv0$ for any continuous process
    of bounded variations on compacts) along with the linearity of the cross variation process,
    and \thmref{1} we have  a.s.~$\P$
   \begin{equation}
   \begin{split}
       &\sd f(S_t)=\sd \left\{f(S_0)+\int_0^tf'(S_s)[dM_s+dV_s]
       +\f12\int_0^tf''(S_s)d\langle M \rangle_s\right\}
      \\
       &= \begin{cases}
       \dis\f12\left[D^+\int_0^tf'(S_s)\,d\langle M,B\rangle_{s}
    +D^-\int_0^tf'(S_s)\,d\langle M,B\rangle_{s}\right]; \ 0<t<\infty,\\
     \dis\left.\left( D^+\int_0^tf'(S_s)\,d\langle M,B\rangle_{s}\right)\right|_{t=0};\ t=0
     \end{cases} \\
     \end{split}
    \label{scrg}
    \end{equation}
    whenever the one sided derivatives are finite.
    If $d\langle M,B\rangle_s=X_{M,B}(s)ds$ and $X_{M,B}$ has continuous paths a.s.~$\P$;
     then it follows from \eqnref{scrg} and the fundamental theorem of classical calculus that
     $\sd f(S_t)=f'(S_t) X_{M,B}(t)$. $t\in\Rp$ a.s.~$\P$.  Also, a.s.~$\P$
     $$X_{M,B}(t)=\f{d}{dt}\int_0^tX_{M,B}(s)ds=\f{d}{dt}\int_0^td\langle M,B\rangle_s=\sd S_t; t\in\Rp,$$
     so that \eqnref{GCR2} follows.
     If $S\in\SBcrep$, then by the same argument above we get that $\sd f(S_t)=f'(S_t)X_t=f'(S_t)\sd S_t$ for all $t\in\Rp$ a.s.~$\P$
     and \eqnref{CR} follows.
     \item From part (i) we have $\sd S^{(2)}(t)=\sd f(B_t)=f'(B_t)$.  Now, using It\^o's formula
     \begin{equation}
     \begin{split}
     \frac{d \langle S^{(2)}\rangle_t}{dt}&=\frac{d}{dt} \left\langle f(B_0)+\int_0^{.} f'(B_s)dB_s+
     \frac12\int_0^{.} f''(B_s)ds\right\rangle_t\\&=\left[f'(B_t)\right]^2
     \end{split}
     \end{equation}
     Also, we have
     \begin{equation}
     \begin{split}
     \frac{d \langle S,S^{(2)}\rangle_t}{dt}&=\frac{d}{dt} \left[\int_0^{t} f'(B_s)d\langle M,B\rangle_{s}\right]
     \\&=\frac{d}{dt}\left[\int_0^{t} f'(B_s)X_{M,B}(s)ds\right]=f'(B_t) X_{M,B}(t)
          \end{split}
     \end{equation}

     So that
     \begin{equation}
     \begin{split}
       \sdS2 S_t=\frac{d \langle S,S^{(2)}\rangle_t}{d \langle S^{(2)}\rangle_t}= \frac{f'(B_t) X_{M,B}(t)}{\left[f'(B_t)\right]^2}
       \end{split}
       \end{equation}
       from which it follows that
       $$\sdS2 S_t\cdot\sd S^{(2)}(t)=\frac{f'(B_t) X_{M,B}(t)}{\left[f'(B_t)\right]^2 }\cdot f'(B_t)=X_{M,B}(t)=
       \frac{d\langle M,B\rangle_t}{dt}=\sd S_t$$
       as claimed.
       \end{enumerate}

      \item If $f:\R\to\R$ is only assumed to be convex, then from standard results in It\^o's calculus
       we get that
       \begin{equation}
f(S_t)=f(S_0)+\int_0^tf'_-(S_s)[dM_s+dV_s]+\f12\Gamma_t^f;\quad t\in\Rp,\mbox{ a.s.}~\P,
 \label{Itoconv}
 \end{equation}
 where $\Gamma^f$ is a continuous increasing process and $f'_-(x)$ is the left derivative of $f$ at $x$.
 The desired results in this convex case then follow by following the same steps in the arguments above, replacing the It\^o formula \eqref{Ito}
 by \eqref{Itoconv} and replacing $f'$ by $f'_-$ throughout.
 \end{enumerate}
 The proof is complete.
       \end{pf*}
       \begin{pf*}{Proof of \thmref{rules}}
       We only need to prove the multiplication and ratio rules \eqref{m} and \eqref{d}, respectively,
       as the addition/subtraction rule is clear by the linearity of the cross variation process.
       We start by proving the multiplication rule \eqref{m}.   To this end, let $f(x,y)=xy$ for $x,y\in\R$.  Applying
       It\^o's formula for functions of several continuous semimartingales (since all the
       $D_if(S^{(1)}_s,S^{(2)}_s)$ and $D_{ij}f(S^{(1)}_s,S^{(2)}_s)$ are continuous) we get:
       \begin{equation}
       \begin{split}
      \sd&\left (S^{(1)}_t S^{(2)}_t\right)=\sd f(S^{(1)}_t,S^{(2)}_t)=\sd f(S^{(1)}_0,S^{(2)}_0)\\
      &+\sd\left[\sum_{i=1}^2 \int_0^tD_if(S^{(1)}_s,S^{(2)}_s)dS_s^{(i)}+\frac12\sum_{0\le i,j\le1}\int_0^tD_{ij}
      f(S^{(1)}_s,S^{(2)}_s)d\langle S^{(i)},S^{(j)}\rangle_s \right]\\
      &=\sd\left[\int_0^tS_s^{(2)}dS_s^{(1)}+\int_0^tS_s^{(1)}dS_s^{(2)}\right]\\
      &=\sd\left[\int_0^tS_s^{(2)}dM_s^{(1)}+\int_0^tS_s^{(1)}dM_s^{(2)}\right]\\
      &=\dis\frac{d}{dt}\left[\int_0^tS_s^{(2)}X_{M^{(1)},B}(s)ds+\int_0^tS_s^{(1)}X_{M^{(2)},B}(s)ds\right]\\
       &= S_t^{(2)}X_{M^{(1)},B}(t)+ S_t^{(1)}X_{M^{(2)},B}(t)
       \\&=S^{(2)}_t\sd S^{(1)}_t+ S^{(1)}_t\sd S^{(2)}_t;\quad0\le t<\infty, \mbox{ a.s.}~\P,
       \end{split}
       \label{mpf}
       \end{equation}
       where we have used \remref{remsd} (a); the assumption on $d\langle M^{(i)},B\rangle_s$ for $i=1,2$;
       the fundamental theorem of deterministic calculus; and the fact that $\sd S^{(i)}_t=X_{M^{(i)},B}(t)$
       for $i=1,2$.

       Now, for the ratio rule \eqref{d} we can either use  the product rule in conjunction with the power rule to show that
       under the additional assumption ($S^{(2)}_t\neq0$ for every $t\in\Rp$, a.s.~$\P$), we have
      \begin{equation*}
      \begin{split}
      \sd \left(S^{(1)}_t\left[S^{(2)}_t\right]^{-1}\right)&=\left[S^{(2)}_t\right]^{-1}\sd S^{(1)}_t-
       S^{(1)}_t\left[S^{(2)}_t\right]^{-2}\sd  S^{(2)}_t\\
       &=\frac{S^{(2)}_t\sd S^{(1)}_t- S^{(1)}_t\sd S^{(2)}_t}{\left[S^{(2)}_t\right]^2};
      \quad0\le t<\infty, \mbox{ a.s.}~\P.
       \end{split}
       \end{equation*}
       Alternatively, we can let $g(x,y)=x/y$ for every $x,y\in\R$ such that $y\neq0$ and proceed exactly as in \eqref{mpf},
       replacing $f(x,y)$ by $g(x,y)$ to get \eqref{d}.
      \end{pf*}

{\bf Acknowledgements.}  This research is supported in part by NSA grant MDA904-03-1-0089.

Department of Mathematical Sciences, Kent State University,
Kent, OH 44242 \\
 Phone: (330) 672-9028 email: allouba@math.kent.edu\\
\end{document}